\documentclass[12pt,a4paper]{article}
\usepackage[latin1]{inputenc}
\usepackage{amsmath,amsthm}
\usepackage{amsfonts}
\usepackage{amssymb}
\usepackage{graphicx,xypic,pstricks}
\input xy
\xyoption{all}

\theoremstyle{definition} 

 \newtheorem{definition}{Definition}[section]

\theoremstyle{plain}      

 \newtheorem{proposition}[definition]{Proposition}
 \newtheorem{theorem}[definition]{Theorem}
 
 \newtheorem{lemma}[definition]{Lemma}

\newcommand{\R}{\mathbb{R}}
\newcommand{\Z}{\mathbb{Z}}
\newcommand{\C}{\mathbb{C}}
\newcommand{\N}{\mathbb{N}}
\newcommand{\Sb}{S'}
\newcommand{\boA}{\mathcal{A}}
\newcommand{\boF}{\mathcal{F}}
\renewcommand{\sl}{\mathrm{Sl}(2,\mathbb{C})}
\DeclareMathOperator{\homo}{Hom}
\DeclareMathOperator{\tr}{Tr}
\DeclareMathOperator{\vect}{span}
\renewcommand{\epsilon}{\varepsilon}
\title{The skein module of torus knots complements}
\author{Julien March\'e
\footnote{Centre de Math\'ematiques Laurent Schwartz,
\'Ecole Polytechnique,
Route de Saclay, 91128 Palaiseau Cedex, France,
email:\,\tt{marche@math.polytechnique.fr}}
}
\date{}

\begin{document}
\maketitle
\begin{abstract}
We compute the Kauffman skein module of the complement of torus knots in $S^3$. Precisely, we show that these modules are isomorphic to the algebra of $\sl$-characters tensored with the ring of Laurent polynomials.
\end{abstract}

Skein modules were introduced indenpendantly by V. Turaev in 1988 and J. Przytycki in 1991 (see \cite{tur1,hp}) as a $\C[A^{\pm 1}]$-module associated to a 3-manifold $M$ generated by banded links inside $M$ with local relations known as Kauffman relations, see for instance \cite{hp}. In the case of $M=S^3$ this construction reduces to the Jones polynomial and in the general case, the evaluation of the skein module at roots of unity is known to fit with the Topological Quantum Field Theory constructed in \cite{bhmv}.
At the same time, skein modules were investigated for themselves. It was shown in \cite{tur2,bul97,prz} that the skein modules of thickened surfaces $\Sigma\times[0,1]$ were non-commutative algebras quantizing in some sense the trace functions on the $\sl$-characters of $\Sigma$. D. Bullock proved in \cite{bul97} that the skein module of $M$ for $A=-1$ was a commutative algebra isomorphic up to nilpotents to the algebra of trace functions on the $\sl$-characters of $M$. Moreover, skein modules were computed in \cite{bl,bul94} for specific 3-manifolds as $S^2\times S^1$, lens spaces, complement of $(2,2p+1)$-torus knots and twist knots. Except for $S^2\times S^1$, skein modules were shown to be free over $\C[A^{\pm 1}]$ and not to have nilpotents when putting $A=-1$. Unfortunately we do not have a general criterium for deciding when this should hold.

Some work of C. Frohman, R. Gelca and S. Garoufalidis \cite{fgl,gar} showed a relation between skein modules and recurrence relations satisfied by colored Jones polynomials. T. Q. T. Lê gave in \cite{le} an efficient description of the skein module of the complement of two-bridge knots. These knots include the previously known examples and the structure of their skein module were shown to play a role in the proof of the AJ-conjecture which relates the family of colored Jones polynomials of a knot $K$ to the $\sl$-character variety of its complement via the $A$-polynomial. This motivates our work on the structure of skein modules of torus knots complements.
We show the following theorem where we denote by $\chi(M)$ the $\sl$-character variety of a manifold $M$, and by $\C[\chi(M)]$ the algebra of regular functions on it.

\begin{theorem}
For all torus knots $T_{p,q}$, the Kauffman skein module of the complement is isomorphic to $\C[A^{\pm 1}]\otimes \C[\chi(S^3\setminus T_{p,q})]$
\end{theorem}
In particular that the skein module is free and do not have nilpotents when $A=-1$. Our proof uses standard tools of skein theory but has some new ingredient: we introduce a filtration on the skein module coming from the intersection number of a link with the separating annulus lying in the complement of torus knots. We find a corresponding filtration in the character algebra and prove our assertion at the graded level.
We hope that this method will give some new informations on skein modules of 3-manifolds containing incompressible surfaces.

After some settings and notations we prove in the first part some lemmas on the relative skein module of a solid torus. This part contains the main technical ingredient of the article. Then we describe the well-known character algebra of torus knots and a degree on it. Finally, we define the filtration on the skein module of torus knots and apply our previous results to compute the associated graded spaces and prove our result.

\section{Settings}
Let $M$ be an oriented 3-manifold with boundary (maybe empty). A banded link in $\partial M$ is by definition a collection $l$ of oriented segments.
A banded link in $M$ relative to $l$ is an embedding of copies of $S^1\times [0,1]$ and $[0,1]\times[0,1]$ into $M$ such that the image of $\{0,1\}\times [0,1]$ coincides with $l$ (with orientation). We will call relative skein module and denote by $S(M,l)$ the free $\C[A^{\pm 1}]$-module generated by isotopy classes of banded links modulo the Kauffman relations (see \cite{prz}). Given a banded link $L$ in $S(M,l)$ homeomorhic to $S^1\times [0,1]$, we will denote by $L^n$ the banded link obtained by replacing $[0,1]$ by $n$ disjoint subintervals.

Let $S^3\subset \C^2$ be defined by the equation $|z|^2+|w|^2=1$. Let $p,q$ be two relatively prime integers:  the torus knot $T_{p,q}$ is the subset of $S^3$ defined by the equation $z^p=w^q$. The sphere $S^3$ can be decomposed into two solid tori $T_1$ and $T_2$ defined respectively by the inequations $|z|^p\le|w|^q$ and $|w|^q\le|z|^p$. These two solid tori intersect on a torus $S$ containing the knot $T_{p,q}$.
Let $X$ be the complement in $S³$ of a small tubular neighborhood of $T_{p,q}$. We will set $X_1=X\cap T_1, X_2=X\cap T_2, R=X\cap S$.

\section{Relative skein modules of solid tori}
Fix a solid torus $T=D^2\times S^1$ and an integer $p\ge 1$. Let $l$ be the curve parametrized by $t\mapsto(e^{it},e^{ipt})$.
For any integer $k\ge 0$, we fix $2k$ real numbers $0\le t_1<t_2<\cdots<t_{2k}<2\pi$ and denote by $S(T,2k)$ the skein module of $T$ relative to the segments $l_i=l(t_i,t_i+\epsilon)$ for $i\in\{1,\ldots,2k\}$ and $\epsilon$ small enough. For simplicity we will refer to $l_i$ as "points" by taking $\epsilon$ arbitrarily small.

We denote by $\Sb(T,2k)$ the quotient of $S(T,2k)$ by the submodule generated by all banded tangles containing an arc parallel to $l$ and joining two end points consecutive in $\R/2\pi\Z$.
We define an automorphism $\tau$ of $S(T,2k)$ in the following way. By an isotopy, move $l_i$ to $l_{i+1}$ for $i<2k$ and $l_{2k}$ to $l_1$ by increasing the parameter $t$.

\begin{lemma}
The automorphism $\tau:S(T,2k)\to S(T,2k)$ satisfies $\tau^{2k}=1$ and induces an automorphism of the quotient $\Sb(T,2k)$.
\end{lemma}
\begin{proof}
Let $\gamma$ be a banded tangle in $S(T,2k)$. Up to isotopy, we can suppose that the intersection of $\gamma$ with the thickened torus $\{|z|>1-\epsilon\}\times S^1$ is a collection of $2k$ segments $s\mapsto ((1-\epsilon s)e^{it_j},(1-\epsilon s)e^{ipt_j})$ for $s\in [0,1]$.
After $2k$ iterations of $\tau$, the end points of $\gamma$ come back to their initial position and the $2k$ segments have been replaced by $s\mapsto ((1-\epsilon s)e^{i(t_j+2\pi s)},(1-\epsilon s)e^{ip(t_j+2\pi s)})$.

There is an isotopy which send this tangle to the initial one. It consists in sliding $T$ along $\gamma$: precisely, the isotopy inside the torus $\{|z|<1-\epsilon\}\times S^1$ is given by the map $\Phi_s(z,w)=(e^{-2i\pi s}z,e^{-2i\pi p s}w)$. In the remaining part of $T$ foliated by tori $\{|z|=\lambda\}$ we interpolate between the map given by $\Phi_s$ on $\lambda=1-\epsilon$ and the identity on $\lambda=1$.
This shows that $\tau^{2k}\gamma$ is isotopic to $\gamma$ and so that $\tau^{2k}$ is the identity.

The second part of the lemma is clear as if $\gamma$ contains an arc joining two consecutive end points and parallel to $l$, then the same will be true for $\tau\gamma$. Hence, $\tau$ factors through the quotient $\Sb(T,2k)$.
\end{proof}

We now give an explicit basis of $S'(T,2k)$ with the description of $\tau$ on it.

We define $e^0_j\in S(T,0)$ in the following way: let $y$ be the banded link $\{0\}\times S^1$ and identify $l$ to a parallel copy of itself pushed into the solid torus. Then, we set $e^0_j=l^n y^m$ were $j=pn+m$ and $m<p$.
\begin{figure}[h]\label{eik}
\begin{center}
\includegraphics[width=5cm]{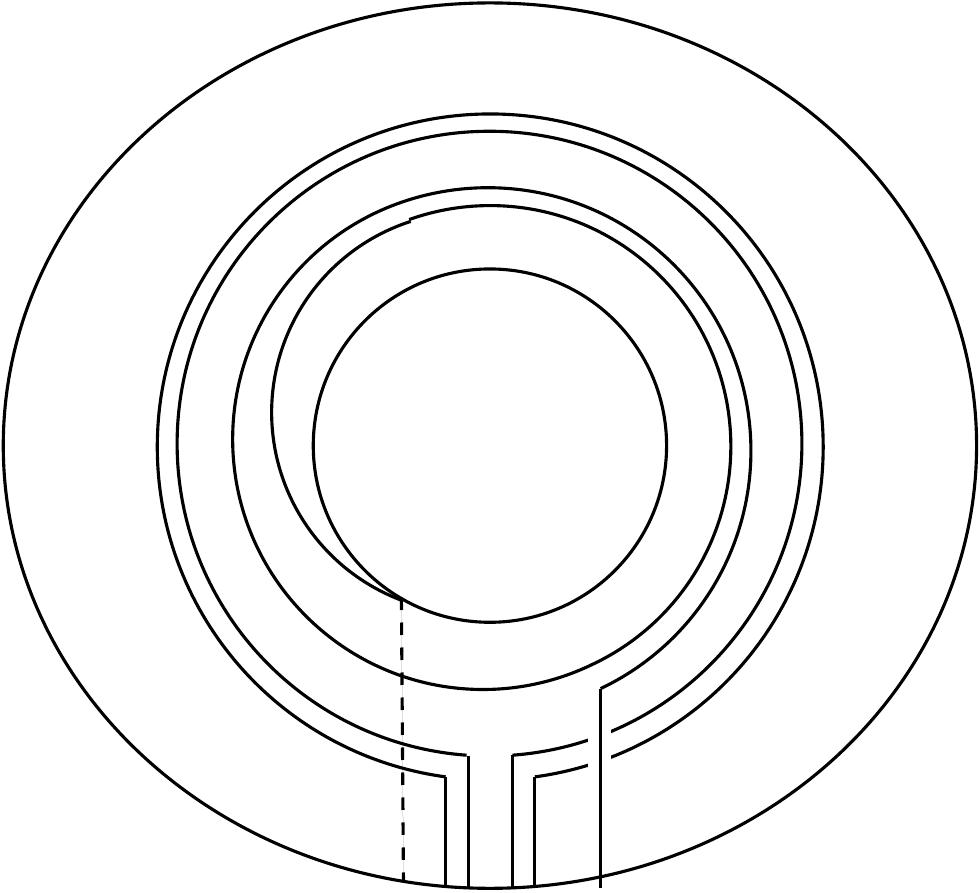}
\caption{The basis element $e^3_2$}
\end{center}
\end{figure}
For $k>0$ and $1\le j\le p-1$, define $e^k_j$ as the union of $k$ components constructed in the following way.
Consider the trivially framed curve $\{1\}\times S^1$, cut a small segment around $(1,1)$, push the remaining part inside $T$ and cable it by $k-1$ parallel strands. We obtain an element in $S(T,2k-2)$ whose end points are numbered $l_2,\ldots,l_{2k-1}$. We add to it the curve $t\mapsto (e^{it},e^{ijt})$ pushed into $T$ with end points $l_1$ and $l_{2k}$. This defines the element $e^k_j$, see Figure \ref{eik} for an example. In the sequel we will modify it by multiplying it with a power of $A$ such that it satisfies the following lemma:

\begin{lemma}\label{lem-fond}
For any $k\in \N$, the elements $e^k_j$ form a basis of $S'(T,2k)$ and for $k>0$, they satisfy the equation $\tau e^k_j=e^k_{p-j}$.
\end{lemma}
\begin{proof}
As a first step, recall that $T$ is homeomorphic to $\Sigma\times [0,1]$ where $\Sigma$ is the annulus $1/2\le|z|\le 1$. We can suppose that the 2k end points are on $S^1\times\{1/2\}$. A standard argument (see \cite{prz}) involving Reidemeister theorem implies that $S(T,2k)$ is the free module generated by isotopy classes of multicurves on $\Sigma$, that is submanifolds of dimension 1 whose boudary coincide with the end points and without closed component bounding a disc.
Using this description, we see that most of the basis elements vanish in $S'(T,2k)$ because all curves in $\Sigma$ joining two end points are either parallel to the boundary or wind once around $\Sigma$. It follows that $S'(T,2k)$ is generated by the elements $y^n$ which by definition are formed by $k$ arcs winding once around $\Sigma$ and $n$ parallel copies of $\{|z|=1/2\}\times \{1/2\}$. We obtain in this way a surjective map $\phi:\Z[A^{\pm 1},y]\to S'(T,2k)$.

Let $z^n$ be the element of $S(T,2k)$ with $k-1$ arcs winding once around $\Sigma$ and a boundary parallel arc between the last two end points. As this element vanishes in $S'(T,2k)$, the same is true for $\tau(z^n)$. Moreover, these vectors generate the kernel of $\phi$ because all the missing relations are presented in this way.

We now invoque another standard result from skein theory: let $\gamma$ be a banded braid with $n$ strands in $[0,1]^3$ joining $X\times \{(1/2,0)\}$ to $X\times \{(1/2,1)\}$ where $X$ is a subset of $[0,1]$ with $n$ elements. Let $I_n$ be the submodule of $S([0,1]^3,2n)$ generated by banded tangles containing boundary parallel arcs. Then there exists $k\in \Z$ such that $\gamma=A^k 1_n \mod I_n$ where $1_n$ is the trivial tangle $X\times \{1/2\}\times[0,1]$.
Applying this argument to the tangle $\tau(z_n)$ which is the closure of a braid with $k-1+n+p$ generators, we find that it decomposes as a polynomial in $z$ of degree $n+p-1$ with leading coefficient invertible in $\Z[A^{\pm 1}]$. Hence the quotient $S'(T,2k)$ is free and generated by $z^j$ for $0\le j<p-1$.

Moreover, the elements $e^k_j$ are by construction closure of braids with $k-1+j$ strands. There is a triangular change of basis with $z^{j-1}$ which shows that $e^k_j$ is also a basis of $S'(T,2k)$.
The same argument shows that $e^0_j$ is a basis of $S(T,0)$ as it is the closure of a braid with $j$ strands.

It remains to compute $\tau e^k_j$ for positive $k$.
Consider the "exterior strand" of $e^k_j$ which winds $j$ times around $T$. If we apply $\tau$, we move its right end point along $l$ in such a way that it winds $p$ times around $T$ and come to the left position. Finally, this strand is a curve joining the two left positions and winding $p-j$ times around $T$. Up to isotopy, the tangles $\tau e^k_j$ and $e^k_{p-j}$ differ only in a neighborhood of the end points. They are closure of different braids with $2k$ strands. The argument above shows that these two elements are proportional up to a power of A, say $\tau e^k_j=A^{u_j}e^k_{p-j}$.

As $\tau^{2k}=1$, we find $(A^{u_j+u_{p-j}})^k=1$ that is $u_{p-j}=-u_j$.
For $j>p/2$, we replace $e^k_j$ with $A^{-u_j}e^k_j$ such that the identity $\tau e^k_{j} =e^k_{p-j}$ holds for all $j$.
\end{proof}

\section{Character variety of torus knots}
The decomposition of $X$ into the subsets $X_1$ and $X_2$ gives a presentation of $\pi_1(X)$ by Van Kampen theorem. The fundamental group of $X_1$ is generated by $u:t\mapsto (e^{it},0)$ and the fundamental group of $X_2$ is generated by $v:t\mapsto (0,e^{it})$. At the intersection, the annulus $R$ has fundamental group generated by $w:t\mapsto(e^{iqt},-e^{ipt})$ assuming that $q$ is even. The path $w$ is homotopic to $u^q$ and $v^p$ so that the fundamental group $G=\pi_1(X)$ is presented by $\langle u,v|u^q=v^p\rangle$.

\subsection{Computation of the character variety}
In this section, we describe the character variety of $G$ denoted by $\chi(G)$, that is the algebraic quotient $\homo(G,\sl)/\!/\sl$. This algebraic variety is well-known (see \cite{kl}) but we give a full description as we need it for our purposes.

It is well known that the character variety of the free group with two generators is isomorphic to $\C^3$ where the isomorphism sends a representation $\rho$ to the triple $(x=\tr\rho(u),y=\tr\rho(v),z=\tr\rho(uv))$. The surjection $\langle u,v\rangle\to G$ gives an injection of $\chi(G)$ into $\C^3$ and we will identify $\chi(G)$ to its image, an algebraic subset of $\C^3$.

The abelianization of $G$ is isomorphic to $\Z$ it is generated by an element $[m]$ such that $[u]=p[m]$ and $[v]=q[m]$. Hence, it is parametrized by the map sending $t$ to $(u(t)=m(t)^p,v(t)=m(t)^q)$ where $m(t)$ is a diagonal matrix with entries $t$ and $1/t$. We obtain $x=t^p+t^{-p},y=t^q+t^{-q},z=t^{p+q}+t^{-p-q}$. Setting $s=t+t^{-1}$, we obtain that abelian characters form an affine curve parametrized by $s\mapsto (T_p(s),T_q(s),T_{p+q}(s))$ where $T_n(t+t^{-1})=t^n+t^{-n}$ are the Tchebychev polynomials.

Let $\rho:G\to\sl$ be an irreducible representation. Then the matrix $\rho(u)^q=\rho(v)^p$ commutes with the image of $\rho$ and hence is equal to $\pm 1$. Hence there is a unique integer $k\in [0,q]$ such that $e^{ki\pi/q}$ is an eigenvalue of $\rho(u)$. As $u$ cannot be equal to $1$ or $-1$, k is in $[1,q-1]$. Again, there is a unique integer $l\in [1,p-1]$ such that $e^{li\pi/p}$ is an eigenvalue of $\rho(v)$.
The equation relating $u$ and $v$ imply that $k$ and $l$ have the same parity. We denote by $\boA$ the set of admissible pairs $(k,l)$: it follows from the discussion that the irreducible part of the character variety $\chi(G)$ is included in $\bigcup_{(k,l)\in \boA}\{\cos(k\pi/q),\cos(l\pi/p)\}\times \C$.

Reciprocally, any element $(x,y,z)$ in this subset correspond to a pair of matrices $(\rho(u),\rho(v))$. As $x=\tr\rho(u)=\cos(k\pi/q)$, we have necessarily $\rho(u)^q=(-1)^k$, and on the other hand $\rho(v)^p=(-1)^l$. We conclude that $(x,y,z)$ corresponds to an element of $\chi(G)$. If this element is abelian, then there is a basis such that $\rho(u)$ and $\rho(v)$ are diagonal with left upper entry equal respectively to $e^{ik\pi/q}$ and $e^{\pm il\pi/p}$. Their product is then diagonal and its trace is equal to $\cos(k\pi/q\pm l\pi/p)$.
This shows that in each line $\{\cos(k\pi/q),\cos(l\pi/p)\}\times \C$, there are two abelian representations and the remaining ones are irreducible.

We can conclude that $\chi(G)$ is a union of $(p-1)(q-1)/2$ disjoint affine lines (irreducible part) and an other line (abelian part) meeting each of the irreducible ones two times. Moreover, we have an explicit description of each of these lines.

\subsection{A filtration on the character algebra}

For $f$ in $\C[\chi(G)]$, we set $\deg f= \max_{(k,l)\in \boA}\deg_z f(\cos(k\pi/q),\cos(l\pi/p),z)$. In other words, $\deg f$ is the maximal degree of the restriction of $f$ to irreducible components.
A function $f$ in $\C[\chi(G)]$ has degree 0 if and only if it is constant on irreducible components. Such a function can be explicitely given by the map  $g(t)=f(t^p+t^{-p},t^q+t^{-q},t^{p+q}+t^{-p-q})$ which is a symetric polynomial in $t$ such that $g(e^{i\pi(k/p+l/q)})=g(e^{i\pi(k/p-l/q)})$ for all $k\in\{1,\ldots,p-1\}$ and $l\in\{1,\ldots,q-1\}$ with the same parity.

On the other hand, the map $\deg: \C[\chi(G)]\to \N$ defines a grading whose graded space for $d>0$ has dimension $(p-1)(q-1)/2$ and is isomorphic to $\C^{\boA}$ via the map sending $f$ to the collection of highest coefficients of $f$ restricted to irreducible components of $\chi(G)$.

We remark that the filtration $F_d$ associated to the degree is more natural as it does not depend on a particular choice of coordinate: this filtration will be given a geometric interpretation in the next section and will be our main ingredient to prove the theorem.

\section{A filtration on the skein module}
Recall that the complement of the torus knot $T_{p,q}$ is decomposed into two parts $X_1$ and $X_2$ glued along an incompressible annulus $R$. Given a banded link $L$ in $X$, one can suppose up to isotopy that it is transverse to $R$. We define its degree as half the geometric intersection number of $L$ with $R$ (which is always an integer as $R$ is separating). Using this notion, we define a filtration of $S(X)$ in the following way:
$$\boF_k(X)=\vect\{ L\in S(X), \deg L\le k\}.$$

\begin{proposition}
There is an isomorphism between $\boF_0$ and $F_0\otimes \C[A^{\pm 1}]$ and between $\boF_k/\boF_{k-1}$ and $(F_k/F_{k-1})\otimes \C[A^{\pm 1}]$ for any $k>0$.
\end{proposition}
Before entering into the proof, we remark that this proposition implies the main theorem of the article. Indeed, if we take any basis of the graded spaces $F_k/F_{k-1}$, the proposition give us a basis of $\boF_k/\boF_{k-1}$ as $\C[A^{\pm 1}]$-module. By lifting these elements to $\boF_k$ for all $k$, we form a basis of $S(X)$ which identifies the skein module with $\C[\chi(X)]\otimes \C[A^{\pm 1}]$ and the proposition is proved.
\begin{proof}
Our strategy is to find elements generating the graduate space of $(\boF_k)$. By setting $A=-1$, we will show that they are independent in $(F_k)$ which will also prove that they are independant in $(\boF_k)$.

We are forced to treat separately the cases $k=0$ and $k>0$.

An element of $\boF_0$ is represented by a disjoint union of banded links in $S(X_1)$ and $S(X_2)$. Using the basis $e^0_j$ of $S(T,0)$, we can present any element of $\boF_0$ as a tensor product $e^0_{j_1}\otimes e^0_{j_2}\in S(X_1)\otimes S(X_2)$. These elements can be written $y_1^{m_1}l_1^{n_1}\otimes
y_2^{m_2}l_2^{n_2}$ with $m_1<p$ and $m_2<q$. But $l_1$ and $l_2$ are parallel copies of $T_{p,q}$ pushed into $X_1$ and $X_2$ which are isotopic to say $l$. Hence, the family $y_1^{m_1}l^n y_2^{m_2}$ generates $\boF_0$ for $n\ge 0, m_1<p$ and $m_2<q$.

Sending $A$ to $-1$, these curves correspond up to sign to the following elements of $\C[G]$:

$\tr(u)^{m_1}\tr(u^q)^n\tr(v)^{m_2}=x^{m_1}P^ny^{m_2}$ where $P=\tr(u^q)=\tr(v^p)$. These functions belong to $F_0$ as they are independant on $z=\tr(uv)$. Moreover, in the variable $t$, they read as $(t^{pm_1}+t^{-pm_1})(t^{pq}+t^{-pq})^n(t^{qm_2}+t^{-qm_2})$. When $m_1,m_2,n$ take all their possible values, the degree of these polynomials take distinct values. This shows that they are linearly independant on $F_0$ and hence on $\boF_0$.

Let $k$ be a positive integer. Any link $L$ in $\boF_k$ can be moved by an isotopy such that the intersection points of $L$ with $R$ lie in a circle parallel to the knot. In other words, the natural map
$\phi:S(X_1,2k)\otimes S(X_2,2k)\to\boF_{k}$
is surjective. Moreover, if the intersection of a link $L$ with one side $X_1$ contains a trivial arc joining two consecutive boundary points, then we can push it by an isotopy into $X_2$ and reduce the degree by 1. This shows that $\phi$ induces an other surjective map

$$\phi:S'(X_1,2k)\otimes S'(X_2,2k)\to\boF_{k}.$$

Next the operator $\tau$ acts on $S'(X_1,2k)$ by a fractional Dehn twist of order $1/2k$. This implies that the link $\phi(\tau x\otimes \tau y)$ is isotopic to $\phi(x\otimes y)$. In other words, we obtain a still surjective map
$$\phi:S'(X_1,2k)\otimes S'(X_2,2k)/\tau \to\boF_{2k}.$$

We now show that this map is an isomorphism by using the explicit basis and computing the associated trace functions.
Let $e^k_{j_1}\otimes e^k_{j_2}$ be the basis of $S'(X_1,2k)\otimes S'(X_2,2k)$ given by Lemma \ref{lem-fond}. By the same lemma, this basis satisfies $\tau(e^k_{j_1}\otimes e^k_{j_2})=e^k_{p-j_1}\otimes e^k_{q-j_2}$.

On the other hand, this basis correspond up to sign to the trace function $\tr(uv)^{k-1}\tr(u^{j_1}v^{j_2})$. We need to compute the corresponding function in $F_{k}/F_{k-1}$, that is we have to compute $\tr(u^iv^j)$ as a polynomial in $x,y,z$ which we do by the following trick using formal power series:
$$G=\sum_{i,j}\tr(u^iv^j)s^it^j=\tr((1-su)^{-1}(1-tv)^{-1})=\frac{\tr((1-su)(1-tv))}{\det(1-su)\det(1-tv)}.$$
The last equality comes from the identity $\tr(M^{-1})=\tr(M)/\det(M)$ for any invertible $2\times 2$ matrix. Then we have $\det(1-su)=(1-s\xi)(1-s\xi^{-1})$ and $\det(1-tv)=(1-t\eta)(1-t\eta^{-1})$ where $\xi$ and $\eta$ are eigenvalues of $u$ and $v$ respectively.

We finally compute: $G=\frac{2-xt-ys+zst}{(1-s\xi)(1-s\xi^{-1})(1-t\eta)(1-t\eta^{-1})}$ from which we deduce that $\tr(u^iv^j)$ is a linear polynomial in $z$ and that its leading coefficient is $\frac{(\xi^i-\xi^{-i})(\eta^j-\eta^{_j})}{(\xi-\xi^{-1})(\eta-\eta^{-1})}$. The proposition results from the following elementary linear algebra lemma:
\begin{lemma}
The map $\C[x^iy^j]_{0<i<p,0<j<q}/\tau\to \C^{\boA}$ sending $x^iy^j$ to $((\xi^i-\xi^{-i})(\eta^j-\eta^{_j}))|_{\xi=e^{ik\pi/p},\eta=e^{il\pi/q}}$ is an isomorphism where $\tau(x^iy^j)=x^{p-i}y^{q-j}$.
\end{lemma}
\begin{proof}
The transformation $\C^{p-1}\to\C^{p-1}$ sending $x^i$ with $0<i<p$ to the vector $(\sin(\pi ij/p))_{0<j<p}$ is well known to be invertible (discrete Fourier transform). Hence, there are unique elements $P_k\in \C[x^i]_{0<i<p}$ and $Q_l\in \C[x^j]_{0<j<q}$ going to the vector $\delta_k\delta_l\in \C^{\boA}$. This shows that the map we are considering is surjective and has to be an isomorphism for dimensional reasons.
\end{proof}
This proves finally the proposition as we have shown that the elements $e^k_i\otimes e^k_j$ generate $\boF_{k}/\boF_{k-1}$ and are linearly independant on $F_k/F_{k-1}$ up to the $\tau$-ambiguity.
\end{proof}

\end{document}